\documentclass[preprint,12pt,sort&compress]{elsarticle}

\usepackage{amsmath}
\usepackage{amsfonts}   
\usepackage{amssymb}  
\usepackage{graphicx}

\newtheorem{thm}{Proposition} 
\newtheorem{lem}[thm]{Lemma}
\newdefinition{rmk}{Remark}
\newproof{pf}{Proof}

\journal{CNSNS}

\begin{document}

\begin{frontmatter}

\title
{IFOHAM-an iterative algorithm based on the first-order equation of HAM: exploratory preliminary results}

\author{Miguel Moreira}

\address
{Portuguese Naval Academy (EN)/ Navy Research Center (CINAV), Base Naval de Lisboa - Alfeite 2810-001 Almada, Portugal}

\begin{abstract}

In this work we present and study an iterative algorithm used to asymptotically solve
nonlinear differential equations. This algorithm (Iterative First Order HAM or
IFOHAM) is based on the first order equation of the Homotopy Analysis Method, HAM.
We show that IFOHAM generalizes Picard-Lindel\"{o}ff 's iteration algorithm.
Moreover, IFOHAM shares with HAM the possibility of ensuring convergence\ by
adequately choosing $c_{0},$ a convergence control parameter. Preliminary
results show that IFOHAM exhibits a very good performance both in aspects
related to the speed of convergence and in aspects related to the CPU
calculation time. It should also be noted that the IFOHAM is a very low
complexity algorithm easily programmable in a symbolic computing
environment.

\end{abstract}


\end{frontmatter}

\section{Introduction}

The Homotopy Analysis Method (HAM) was introduced in 1992 by Shijun Liao in
its PhD thesis \cite{Liao0} and subsequently developed and applied by this
author \cite{Liao1,liao2,liao3} and by a growing community of researchers.

An extensive and complete state of the art concerning the HAM technique can be found in
\cite{liao3}.

This technique inserts or relates to the so-called asymptotic methods
\cite{bayat1} and analytical approximation methods \cite{rhadika1,zwillinger1}.

Basically, the HAM technique transforms the original nonlinear problem
(nonlinear differential equation, nonlinear differential equation system or
even nonlinear partial differential equation system, for instance) into a set
of linear differential equations (to be solved recursively) whose analytic
solutions constitute the terms of a series of functions representing the
solution of the original problem. This transformation is based on the concept
of homotopy (under which an initial guess of the solution is continuously
deformed to the solution of the original equation) and is built from the
so-called zeroth-order deformation equation.

Consider the Initial Value Problem (IVP) described by
\begin{equation}
{N}\left[  x\left(  t\right)  \right]  =0 \label{eq_0010},
\end{equation}
where $N$ represents a nonlinear operator. So, depending on the order of the
problem, the solution $x=x\left(  t\right)  $ must satisfy initial conditions,
such as
\begin{equation}
\left \{
\begin{array}
[c]{c}
x\left(  t_{0}\right)  =x_{0}^{\left(  0\right)  }\\
x^{\left(  1\right)  }\left(  t_{0}\right)  =x_{0}^{\left(  1\right)  }\\
\vdots \\
x^{\left(  r-1\right)  }\left(  t_{0}\right)  =x_{0}^{\left(  r-1\right)  }
\end{array}
\right.  , \label{eq_0011}
\end{equation}
if we assume that (\ref{eq_0010}) is defined by an ordinary differential
equation of order $r$.

This work will be centered on the basic formulation of HAM developed in \cite{liao4} which is supported %
by the corresponding zeroth-order deformation equation
(\ref{eq_0035}). Based on the previously mentioned equation an iterative algorithm (iterative
first-order HAM: IFOHAM) to solve (\ref{eq_0010}) will be proposed and its
main features will be presented and discussed.

\section{Basic HAM}

\subsection{Zeroth-order deformation equation}

Following \cite{liao4}, a zeroth-order deformation equation (\ref{eq_0035}) is
defined, where $\mathcal{L}$ %
represents an appropriate linear operator, $c_{0}\neq0$ stands for the
convergence control parameter of HAM (to be described later), $q\in \left[
0,1\right]  $ and $N$ represents the non-linear operator describing the
problem (\ref{eq_0010}) to be solved:
\begin{equation}
\left(  1-q\right)%
\mathcal{L}
\left[  \phi \left(  t;q\right)  -u_{0}\left(  t\right)  \right]
=c_{0}q\left \{  N\left[  \phi \left(  t;q\right)  \right]
\right \}  .\label{eq_0035}
\end{equation}
In expression (\ref{eq_0035}), $\phi=\phi \left(  t;q\right)  $ represents the so-called homotopy
MacLaurin series which is a power series of \ the embedding parameter
$q$:
\begin{equation}
\phi \left(  t;q\right)  =u_{0}\left(  t\right)  +\sum_{n=1}^{+\infty}
u_{n}\left(  t\right)  q^{n},\text{ }q\in \left[  0,1\right]  .\label{eq_0040}
\end{equation}
Observe that in (\ref{eq_0035}), $u_{0}\left(  t\right)  $ represents an initial
guess (to be postulated), satisfying the initial conditions of the solution,
$u=u\left(  t\right)  $, of our original problem (\ref{eq_0010}). Note that
$u_{0}\left(  t\right)  $ is also the zeroth-order term of the homotopy
Maclaurin series (\ref{eq_0040}), that is,
\begin{equation}
\phi \left(  t;0\right)  =u_{0}\left(  t\right)  .\label{eq_0045}
\end{equation}
Setting $q=0,$ in the zeroth-order deformation equation (\ref{eq_0035}), we
obtain
$$\mathcal{L}
\left[  \phi \left(  t;0\right)  -u_{0}\left(  t\right)  \right]  =
\mathcal{L}
\left[  0\right]  =0.$$
Setting $q=1$ we obtain $N\left[  \phi \left(  t;1\right)
\right]  =0$. This fact shows that, converging
\begin{equation}
\phi \left(  t;1\right)  =u_{0}\left(  t\right)  +\sum_{n=1}^{+\infty}
u_{n}\left(  t\right)  ,\label{eq_0050}%
\end{equation}
(\ref{eq_0050}) is solution of (\ref{eq_0010}). So, the coefficients of the
homotopy MacLaurin series (\ref{eq_0040}) are precisely the terms
$u_{n}\left(  t\right)  ,$ $n\in
\mathbb{N}
_{0}$, of the series of functions representing the searched solution
\begin{equation}
u\left(  t\right)  =u_{0}\left(  t\right)  +\sum_{n=1}^{+\infty}u_{n}\left(
t\right) \label{eq_0055}
\end{equation}
of our problem (\ref{eq_0010}).

Typically the zeroth-order deformation equation (\ref{eq_0035}) is indexed in
the parameter $q$ (embedding parameter) and constitutes an homotopic family of
differential equations with homotopic solutions under the embedding parameter
$q,$ each one described by $\phi=\phi \left(  t;q\right)  $. If
$q=0,$ then $\phi \left(  t;0\right)  =u_{0}\left(  t\right)  ,$ will be the
trivial solution of
\begin{equation}
\mathcal{L}
\left[  u\left(  t\right)  \right]  =
\mathcal{L}
\left[  u_{0}\left(  t\right)  \right]  \text{.}
\end{equation}
If $q=1$, then
$
N\left[  \phi \left(  t;1\right)  \right]  =0
$, \text{ that is, } $\phi \left(  t;1\right)  $ will be our searched solution.

It should be noted that the zeroth-order deformation equation (\ref{eq_0035}) is the
starting point of HAM. Generalized formulations of the zeroth-order
deformation equation can be built to be applied more efficiently as well as in
broader contexts, \cite{liao2,liao3}.

HAM users are interested in $\phi \left(  t;1\right)  $, that is, in the
solution of problem (\ref{eq_0010}). Let's see how (\ref{eq_0050}) can be
obtained using this method.

\subsection{High order deformation equations}

Define the operator%
\begin{equation}
\mathcal{D}_{k}=\left.  \frac{1}{k!}\frac{\partial^{k}}{\partial q^{k}%
}\right \vert _{q=0}\label{eq_00100}%
\end{equation}
and let's apply it to the zeroth-order deformation (\ref{eq_0035}). One obtain
(see \cite{liao2}):%
\begin{equation}%
\mathcal{L}%
\left[  u_{1}\left(  t\right)  \right]  =c_{0}\left[  N\left[
u_{0}\left(  t\right)  \right]  \right] \label{eq_def_ordem_1}%
\end{equation}
and%
\begin{equation}%
\mathcal{L}%
\left[  u_{n}\left(  t\right)  -u_{n-1}\left(  t\right)  \right]
=c_{0}\mathcal{D}_{n-1}\left[  N\left[  \phi \left(  t;q\right)
\right]  \right]  \text{, }n\in%
\mathbb{N}
\text{ and }n\geq2.\label{eq_def_ordem_superior}%
\end{equation}

Equations (\ref{eq_def_ordem_1}) and (\ref{eq_def_ordem_superior}) constitutes
the\ so-called high order deformation equations. These equations are linear
and can be recursively solved to obtain each term $u_{n}\left(  t\right)  ,$
$n\in%
\mathbb{N}
_{0}$ of (\ref{eq_0050}). Typically, using a symbolic computer environment,
such as Mathematica, Maple or Matlab, for instance, one can automatically
solve (\ref{eq_def_ordem_1}) and (\ref{eq_def_ordem_superior}) and obtain an
approximate solution
\begin{equation}
u^{m}\left(  t\right)  =u_{0}\left(  t\right)  +\sum_{n=1}^{m}u_{n}\left(
t\right) \label{eq_00110}%
\end{equation}
of order $m$ of the problem (\ref{eq_0010}). This approximate solution can be
called $m$th-order solution.

Note additionally that (\ref{eq_00110}) must satisfy the initial conditions
(\ref{eq_0011}) of our problem and $u_{0}\left(  t\right)  $ already does.
Therefore $u_{n}\left(  t\right)  $ and their derivatives up to order $r-1$
must satisfy null initial conditions for $n=1,\ldots,m$. In summary:%
\begin{equation}
\left \{
\begin{array}
[c]{c}%
u_{0}\left(  t_{0}\right)  =x_{0}^{\left(  0\right)  }\\
u_{0}^{\left(  1\right)  }\left(  t_{0}\right)  =x_{0}^{\left(  1\right)  }\\
\vdots \\
u_{0}^{\left(  r-1\right)  }\left(  t_{0}\right)  =x_{0}^{\left(  r-1\right)
}%
\end{array}
\right.  \text{ and }\left \{
\begin{array}
[c]{c}%
u_{n}\left(  t_{0}\right)  =0\\
u_{n}^{\left(  1\right)  }\left(  t_{0}\right)  =0\\
\vdots \\
u_{n}^{\left(  r-1\right)  }\left(  t_{0}\right)  =0
\end{array}
\right.  ,\forall n=1,\ldots,m\label{eq_00120}%
\end{equation}

\subsection{A\ trivial example of application of HAM}

For the sake of clarity in exposition let's apply the HAM technique to a
nonlinear initial value problem with the known the solution $x=\tan t$:%
\begin{equation}
\left \{
\begin{array}
[c]{c}%
x^{\prime}=1+x^{2}\\
x\left(  0\right)  =0
\end{array}
\right.  .\label{eq_001010}%
\end{equation}

This IVP will be also used later as a simple test case.

Let's consider%
\begin{equation}
N[x]=x^{\prime}-x^{2}-1,\label{eq_001005}%
\end{equation}
define
\[%
\mathcal{L}%
\left[  \cdot \right]  =\frac{d}{dt}\left[  \cdot \right]  ,
\]
consider the convergence control parameter $c_{0}=-1,$ define the homotopy
Maclaurin series
\[
\phi \left(  t;q\right)  =u_{0}\left(  t\right)  +\sum_{n=1}^{+\infty}%
u_{n}\left(  t\right)  q^{n},\text{ }q\in \left[  0,1\right]  ,
\]
and choose the following initial guess (satisfying the initial conditions)%
\begin{equation}
u_{0}\left(  t\right)  =t.\label{eq_001040}%
\end{equation}
Hence, the zeroth-order deformation equation is%
\begin{equation}
\left(  1-q\right)  \frac{d}{dt}\left[  \phi \left(  t;q\right)  -t\right]
=-q\left \{  \frac{\partial \phi \left(  t;q\right)  }{\partial t}-\left[
\phi \left(  t;q\right)  \right]  ^{2}-1\right \}  \text{, }q\in \left[
0,1\right]  .\label{eq_001050}%
\end{equation}
and the corresponding high-order homotopy equations are%
\begin{gather}
\frac{d}{dt}\left[  u_{m}\left(  t\right)  -\chi_{m}u_{m-1}\left(  t\right)
\right]  =-\mathcal{D}_{m-1}\left[  \frac{\partial \phi \left(  t;q\right)
}{\partial t}-\left[  \phi \left(  t;q\right)  \right]  ^{2}-1\right]
,\label{eq_001060}\\
m\geq1\text{ and }\chi_{m}=\left \{
\begin{array}
[c]{c}%
0\text{ if }m=1\\
1\text{ if }m>1
\end{array}
\right.  .\nonumber
\end{gather}
Applying (\ref{eq_00100}) one deduce from (\ref{eq_001060}) the following
high-order deformation equations:%
\begin{gather}
\frac{d}{dt}\left[  u_{m}\left(  t\right)  \right]  =\left(  \chi
_{m}-1\right)  \left(  \frac{d}{dt}\left[  u_{m-1}\left(  t\right)  \right]
-1\right)  +%
{\displaystyle \sum \limits_{k=0}^{m-1}}
u_{k}\left(  t\right)  u_{m-1-k}\left(  t\right)  ,\label{eq_001070}\\
m\geq1\text{ and }\chi_{m}=\left \{
\begin{array}
[c]{c}%
0\text{ if }m=1\\
1\text{ if }m>1
\end{array}
\right.  .\nonumber
\end{gather}
From (\ref{eq_001040}) and (\ref{eq_001070}) let's present the first four
linear ordinary differential equations as well as the corresponding solutions
recursively solved:%
\begin{equation}%
\begin{tabular}
[c]{l|l|l}\hline \hline
$m$ & Linear equation & Solution\\ \hline \hline
$1$ & $\frac{du_{1}}{dt}=t^{2}$ & $u_{1}=\frac{t^{3}}{3}$\\
$2$ & $\frac{du_{2}}{dt}=2tu_{1}$ & $u_{2}=\frac{2t^{5}}{15}$\\
$3$ & $\frac{du_{3}}{dt}=u_{1}^{2}+2tu_{2}$ & $u_{3}=\frac{17t^{7}}{315}$\\
$4$ & $\frac{du_{4}}{dt}=2u_{1}u_{2}+2tu_{3}$ & $u_{4}=\frac{62t^{9}}{2835}$%
\end{tabular}
\label{eq_001080}%
\end{equation}
Based on (\ref{eq_001080}) one can write the fourth order solution of the IVP
(\ref{eq_001010}):%
\begin{equation}
u^{4}\left(  t\right)  =t\allowbreak+\frac{1}{3}t^{3}+\frac{2}{15}t^{5}%
+\frac{17}{315}t^{7}+\frac{62}{2835}t^{9}.\label{eq_001090}%
\end{equation}
Observe and compare (\ref{eq_001090}) with the Maclaurin series of $\tan t:$%
\begin{align*}
\tan t  & =t\allowbreak+\frac{1}{3}t^{3}+\frac{2}{15}t^{5}+\frac{17}{315}%
t^{7}+\frac{62}{2835}t^{9}\\
& +\frac{1382}{155\,925}\allowbreak t^{11}+\frac{21\,844}{6081\,075}%
t^{13}+\frac{929\,569}{638\,512\,875}t^{15}+O\left(  t^{17}\right)  .
\end{align*}
It can be stated that HAM \textquotedblleft surgically \textquotedblright \ determines
the terms of the Maclaurin series of the solution of our problem.

\subsection{Main features of HAM}

All the information needed to find the terms of (\ref{eq_0050}) are contained
in the zeroth-order equation (\ref{eq_0035}). One important parameter in this
equation, see \cite{Liao0,Liao1,liao2,liao3}, is
precisely $c_{0}$ which controls the convergence/divergence of the series
solution of (\ref{eq_0010}). This parameter is called convergence control
parameter and need to be carefully chosen. In \cite{Liao0,Liao1,liao2} %
are presented some practical approaches to choose $c_{0}$ in
order to ensure the convergence as well as the speed of convergence of the
series solution built in the frame of HAM. Besides, the user of HAM has a
great freedom in choosing the linear operator $%
\mathcal{L}%
$ as well as the initial guess, $u_{0}\left(  t\right)  ,$ of the solution.
All these facts underlies some remarkable advantages of HAM, namely:

\begin{enumerate}
\item Guarantee of convergence by adequately choosing $c_{0},$ the convergence
control parameter;

\item Flexibility on the choice of base functions and decide about the
solution expression by adequately choosing $%
\mathcal{L}%
$ and the initial guess $u_{0}\left(  t\right)  $;
\item Ability to find the main parameters, such as amplitude and frequency, of
periodic solutions of nonlinear evolution problems;

\item Great generality of application ranging from solving weakly to strong
nonlinear differential equations or even fractional differential equations.

\end{enumerate}
\section{IFOHAM-Iterative first order HAM}

\subsection{Motivation}

Consider the original IVP problem (\ref{eq_0010}) and (\ref{eq_0011}). Suppose
that (\ref{eq_0050}) converges and consider the first-order deformation
equation (\ref{eq_def_ordem_1})%
\[%
\mathcal{L}%
\left[  u_{1}\left(  t\right)  \right]  =c_{0}\left[  N\left[
u_{0}\left(  t\right)  \right]  \right]  ,
\]
from which we can obtain $u_{1}$. It will be reasonable to conjecture that
$u_{0}\left(  t\right)  +u_{1}\left(  t\right)  $ will be a best initial guess
than the (postulated) original one $u_{0}\left(  t\right)  $. This argument
suggest the following iterative procedure to improve the initial guess $u_{0}
$ for the solution of (\ref{eq_0010}):%
\begin{equation}%
\mathcal{L}%
\left[  u_{m+1}\left(  t\right)  \right]  =c_{0}\left[  N\left[
\sum_{k=0}^{m}u_{k}\left(  t\right)  \right]  \right]  \text{, }%
m\geq0.\label{eq_itfoham_v1}%
\end{equation}

As was the case in applying HAM, in accordance with (\ref{eq_00120}), one must
assure that $u_{k}\left(  t\right)  $ and their derivatives up to order $r-1$
must satisfy null initial conditions for $k=1,\ldots,m+1$. For instance, if
$N[\cdot]$ is defined by a first-order nonlinear differential equation, then%
\begin{equation}
\left \{
\begin{array}
[c]{c}%
u_{0}\left(  t_{0}\right)  =x_{0}^{\left(  0\right)  }\\
u_{1}\left(  t_{0}\right)  =0\\
\vdots \\
u_{m}\left(  t_{0}\right)  =0\\
\vdots
\end{array}
\right.  .\label{eq_002004}%
\end{equation}

Algorithm (\ref{eq_itfoham_v1}) is entirely based on the first-order
deformation equation (\ref{eq_def_ordem_1}) of HAM. So, let's call it
iterative first-order HAM: IFOHAM.

Define
\begin{equation}
x_{m}\left(  t\right)  =\sum_{k=0}^{m}u_{k}\left(  t\right)
.\label{eq_itfoham_v2}%
\end{equation}
and call (\ref{eq_itfoham_v2}) an $m$th-order solution of the problem
[(\ref{eq_0010}) and (\ref{eq_0011})].

Some interesting issues arise immediately:

\begin{enumerate}
\item Does (\ref{eq_itfoham_v2}) converge for the solution of (\ref{eq_0010})?
In what circunstances?

\item How does compare or relate (\ref{eq_itfoham_v1}) with other iterative algorithms?

\item How does the performance of (\ref{eq_itfoham_v1}) relates to the
performance of HAM?

\item What features (\ref{eq_itfoham_v1}) share with HAM? In what features is
(\ref{eq_itfoham_v1}) better effective than HAM?
\end{enumerate}

In the following we will respond these issues and we will present some
exploratory preliminary results.

\subsection{IFOHAM and Picard-Lindel\"{o}ff 's iteration algorithm}

Consider the IVP described in the following first-order ordinary differential
equation and the corresponding initial condition:
\begin{equation}
\left \{
\begin{array}
[c]{c}%
\frac{dx}{dt}=f\left(  t,x\right) \\
x\left(  t_{0}\right)  =x_{0}^{\left(  0\right)  }%
\end{array}
\right.  .\label{eq_002005}%
\end{equation}
Note that in this case the nonlinear operator $N[\cdot]$ can be identified
with an ordinary differential equation in the canonical form, that is
\begin{equation}
N\left[  x\right]  \equiv \frac{dx}{dt}-f\left(  t,x\right)
.\label{eq_002010}%
\end{equation}
Due to (\ref{eq_002010}), IFOHAM (\ref{eq_itfoham_v1}) reduces to%
\begin{equation}%
\mathcal{L}%
\left[  u_{m+1}\left(  t\right)  \right]  =c_{0}\left[  \sum_{k=0}^{m}%
u_{k}^{\prime}\left(  t\right)  -f\left(  t,\sum_{k=0}^{m}u_{k}\left(
t\right)  \right)  \right]  \text{, }m\geq0.\label{eq_002030}%
\end{equation}

Let
\begin{equation}
u_{0}\left(  t\right)  =x_{0}^{\left(  0\right)  },\label{eq_002031}%
\end{equation}
be our initial guess, and assume
\begin{equation}
u_{k}\left(  t_{0}\right)  =0\text{, }\forall k\in%
\mathbb{N}
.\label{eq_002032}%
\end{equation}
Consider
\begin{equation}
c_{0}=-1,\label{eq_002035}%
\end{equation}
and define
\begin{equation}%
\mathcal{L}%
\left[  h\left(  t\right)  \right]  =\frac{dh}{dt}\left(  t\right)  \text{ and
}%
\mathcal{L}%
^{-1}\left[  h\left(  t\right)  \right]  =\int_{t_{0}}^{t}h\left(  \xi \right)
d\xi.\label{eq_002040}%
\end{equation}
From (\ref{eq_002030}) using (\ref{eq_002004}), (\ref{eq_002031}),
(\ref{eq_002035}) and (\ref{eq_002040}) one deduce
\begin{equation}
\sum_{k=0}^{m+1}u_{k}\left(  t\right)  =x_{0}^{\left(  0\right)  }+\int
_{t_{0}}^{t}f\left(  \xi,\sum_{k=0}^{m}u_{k}\left(  \xi \right)  \right)
d\xi,\label{eq_002050}%
\end{equation}
That is,%
\begin{equation}
\left \{
\begin{array}
[c]{l}%
x_{0}\left(  t\right)  =x_{0}^{\left(  0\right)  }\\
x_{m+1}\left(  t\right)  =x_{0}^{\left(  0\right)  }+\int_{t_{0}}^{t}f\left(
\xi,x_{m}\left(  \xi \right)  \right)  d\xi,\text{ }m\geq0
\end{array}
\right.  ,\label{eq_002055}%
\end{equation}
where
\begin{equation}
x_{m}\left(  t\right)  =\sum_{k=0}^{m}u_{k}\left(  t\right)
.\label{eq_002070}%
\end{equation}
Clearly, (\ref{eq_002055}) represents Picard-Lindel\"{o}ff 's iterative
algorithm. So, in this case and under the described restritions IFOHAM
(\ref{eq_002030}) and Picard-Lindel\"{o}ff 's iteration algorithm
(\ref{eq_002055}) generate the same sequence of functions. The following
result can be stated:

\begin{thm}
\label{proposition1}Consider the IVP
\[
\left \{
\begin{array}
[c]{c}%
\frac{dx}{dt}=f\left(  t,x\right) \\
x\left(  t_{0}\right)  =x_{0}^{\left(  0\right)  }%
\end{array}
\right.  ,
\]
where $f$ is a continuous real function on an open set $D\in%
\mathbb{R}
^{2}$ and suppose that $\left(  t_{0},x_{0}^{\left(  0\right)  }\right)  \in
D$. Consider additionally the corresponding nonlinear operator%
\[
N\left[  x\right]  \equiv \frac{dx}{dt}-f\left(  t,x\right)  ,
\]
the IFOHAM algorithm
\begin{equation}
\left \{
\begin{array}
[c]{l}%
u_{0}\left(  t\right)  =x_{0}^{0}\\%
\mathcal{L}%
\left[  u_{m+1}\left(  t\right)  \right]  =c_{0}\left[  N\left[
\sum_{k=0}^{m}u_{k}\left(  t\right)  \right]  \right]  ,\text{ }m\geq0\\
\text{with }u_{k}\left(  t_{0}\right)  =0\text{, }\forall k\in%
\mathbb{N}%
\end{array}
\right.  ,\label{eq_p_0005}%
\end{equation}
where $c_{0}=-1,$ $%
\mathcal{L}%
\left[  h\left(  t\right)  \right]  =\frac{dh}{dt}\left(  t\right)  $ and $%
\mathcal{L}%
^{-1}\left[  h\left(  t\right)  \right]  =\int_{t_{0}}^{t}h\left(  \xi \right)
d\xi$ and consider further the Picard-Lindel\"{o}ff 's iteration algorithm
\begin{equation}
\left \{
\begin{array}
[c]{l}%
x_{0}\left(  t\right)  =x_{0}^{\left(  0\right)  }\\
x_{m+1}\left(  t\right)  =x_{0}^{\left(  0\right)  }+\int_{t_{0}}^{t}f\left(
\xi,x_{m}\left(  \xi \right)  \right)  d\xi,\text{ }m\geq0
\end{array}
\right.  .\label{eq_p_0020}%
\end{equation}
\ 

Then,
\begin{equation}
x_{m}\left(  t\right)  =\sum_{k=0}^{m}u_{k}\left(  t\right)  \text{, }\forall
m\in%
\mathbb{N}
_{0}\label{eq_p_0030}%
\end{equation}
whenever $\left(  t,x_{k}\left(  t\right)  \right)  \in D$ for $k=1,\ldots
,m-1.$

\begin{pf}
This statement is the instance $c_{0}=-1$ of Proposition \ref{proposition_4}%
.\hfill$\square$
\end{pf}
\end{thm}

This means that under the described restrictions the convergence of IFOHAM\ is
ensured if (\ref{eq_002005}) satisfies the classical Picard-Lindel\"{o}ff 's
conditions for the existence and uniqueness of a solution. In short:

\begin{thm}
\label{proposition2}Let $D$ be an open set in $%
\mathbb{R}
^{2}$. Let $\left(  t_{0},x_{0}^{\left(  0\right)  }\right)  \in D$ and let
$a$ and $b$ be positive constants such that the set%
\[
R=\left \{  \left(  t,x\right)  :\left \vert t-t_{0}\right \vert \leq a\text{ and
}\left \vert x-x_{0}^{\left(  0\right)  }\right \vert \leq b\right \}
\]
is contained in $D.$ Suppose function $f$ is continuous and defined on $D$ and
satisfies a Lipschitz condition with respect to $x$ in $R$. Let $M=\underset
{\left(  t,x\right)  \in R}{\max}\left \vert f\left(  t,x\right)  \right \vert $
and $A=\min \left \{  a,\frac{b}{M}\right \}  .$ Then,

\begin{description}
\item[a)] the the IVP
\[
\left \{
\begin{array}
[c]{c}%
\frac{dx}{dt}=f\left(  t,x\right) \\
x\left(  t_{0}\right)  =x_{0}^{\left(  0\right)  }%
\end{array}
\right.  ,
\]
has a unique solution $x=x\left(  t\right)  $ on the open interval $I=\left]
t_{0}-A,t_{0}+A\right[  $.

\item[b)] the sequence $\left \{  x_{n}\left(  t\right)  \right \}  $, where
\[
\left \{
\begin{array}
[c]{l}%
x_{0}\left(  t\right)  =x_{0}^{\left(  0\right)  }\\
x_{n+1}\left(  t\right)  =x_{0}^{\left(  0\right)  }+\int_{t_{0}}^{t}f\left(
\xi,x_{n}\left(  \xi \right)  \right)  d\xi,\text{ }n\geq0
\end{array}
\right.  ,
\]
converges uniformly on $I$ to $x\left(  t\right)  $;

\item[c)] the sequence $\left \{  \sum_{k=0}^{n}u_{k}\left(  t\right)
\right \}  $ where
\[
\left \{
\begin{array}
[c]{l}%
u_{0}\left(  t\right)  =x_{0}^{0}\\%
\mathcal{L}%
\left[  u_{n+1}\left(  t\right)  \right]  =c_{0}\left[  N\left[
\sum_{k=0}^{n}u_{k}\left(  t\right)  \right]  \right]  ,\text{ }n\geq0\\
\text{with }u_{k}\left(  t_{0}\right)  =0\text{, }\forall k\in%
\mathbb{N}%
\end{array}
\right.  ,
\]
whith $N\left[  x\right]  \equiv \frac{dx}{dt}-f\left(
t,x\right)  $, $c_{0}=-1,$ $%
\mathcal{L}%
\left[  h\left(  t\right)  \right]  =\frac{dh}{dt}\left(  t\right)  $ and $%
\mathcal{L}%
^{-1}\left[  h\left(  t\right)  \right]  =\int_{t_{0}}^{t}h\left(  \xi \right)
d\xi$ converges uniformly on $I $ to $x\left(  t\right)  $.
\end{description}

\begin{pf}
The proof of parts (a) and (b) are classic and can be found for instance in
\cite{Cronin1}. Part (c) is an immediate consequence of Proposition
\ref{proposition1}.\hfill$\square$
\end{pf}
\end{thm}

Let us now study the role of the convergence control parameter $c_{0}$ in the
behavior of IFOHAM.

\subsection{IFOHAM and the convergence control parameter $c_{0}$}

Consider again the IVP described in (\ref{eq_002005}) and let the
corresponding nonlinear operator $N[\cdot]$ be
\begin{equation}
N\left[  x\right]  \equiv \frac{dx}{dt}-f\left(  t,x\right)
.\label{eq_300005}%
\end{equation}
So, IFOHAM (\ref{eq_itfoham_v1}) reduces to%
\begin{equation}%
\mathcal{L}%
\left[  u_{m+1}\left(  t\right)  \right]  =c_{0}\left[  \sum_{k=0}^{m}%
u_{k}^{\prime}\left(  t\right)  -f\left(  t,\sum_{k=0}^{m}u_{k}\left(
t\right)  \right)  \right]  \text{, }m\geq0.\label{eq_300010}%
\end{equation}
Let our initial guess be
\begin{equation}
u_{0}\left(  t\right)  =x_{0}^{\left(  0\right)  },\label{eq_300020}%
\end{equation}
and define
\begin{equation}%
\mathcal{L}%
\left[  h\left(  t\right)  \right]  =\frac{dh}{dt}\left(  t\right)  \text{ and
}%
\mathcal{L}%
^{-1}\left[  h\left(  t\right)  \right]  =\int_{t_{0}}^{t}h\left(  \xi \right)
d\xi.\label{eq_300030}%
\end{equation}

From (\ref{eq_300010}) one deduce using (\ref{eq_002004}), (\ref{eq_300020})
and (\ref{eq_300030}),%
\begin{equation}
\sum_{k=0}^{m+1}u_{k}\left(  t\right)  =\left(  1+c_{0}\right)  \sum_{k=0}%
^{m}u_{k}\left(  t\right)  -c_{0}\left[  x_{0}^{\left(  0\right)  }%
+\int_{t_{0}}^{t}f\left(  \xi,\sum_{k=0}^{m}u_{k}\left(  \xi \right)  \right)
d\xi \right]  ,\label{eq_300040}%
\end{equation}
or, equivalently using (\ref{eq_002070})%
\begin{equation}
x_{m+1}\left(  t\right)  =\left(  1+c_{0}\right)  x_{m}-c_{0}\left[
x_{0}^{\left(  0\right)  }+\int_{t_{0}}^{t}f\left(  \xi,x_{m}\left(
\xi \right)  \right)  d\xi \right]  .\label{eq_300050}%
\end{equation}

Note that interestingly (\ref{eq_300050}) can be interpreted as a weighted
average between $x_{m}$, the previous iteration, and
\[
x_{0}^{\left(  0\right)  }+\int_{t_{0}}^{t}f\left(  \xi,x_{m}\left(
\xi \right)  \right)  d\xi,
\]
the iterate $m+1$ computed using Picard-Lindel\"{o}ff 's iterative algorithm
(\ref{eq_002055}). This fact suggest the decrease of the convergence speed of
the algorithm for increasing values of $c_{0}$ in the interval $\left[
-1,0\right[  $. In reality this conjecture will be corroborated by expression
(\ref{eq_p5_002020}) from Proposition \ref{proposition_5}.

From (\ref{eq_300050}) one can deduce the equivalent useful expression:%
\begin{equation}
x_{m+1}=x_{0}^{\left(  0\right)  }+\left(  1+c_{0}\right)  \int_{t_{0}}%
^{t}\frac{dx_{m}}{dt}\left(  \xi \right)  d\xi-c_{0}\int_{t_{0}}^{t}f\left(
\xi,x_{m}\left(  \xi \right)  \right)  d\xi.\label{eq_300060}%
\end{equation}
Let's summarize these results:

\begin{thm}
\label{proposition_4}Consider the IVP
\[
\left \{
\begin{array}
[c]{c}%
\frac{dx}{dt}=f\left(  t,x\right) \\
x\left(  t_{0}\right)  =x_{0}^{\left(  0\right)  }%
\end{array}
\right.  ,
\]
where $f$ is a continuous real function on an open set $D\in%
\mathbb{R}
^{2}$ and suppose that $\left(  t_{0},x_{0}^{\left(  0\right)  }\right)  \in
D$. Let $c_{0}\in%
\mathbb{R}
$ and consider algorithm
\begin{equation}
\left \{
\begin{array}
[c]{l}%
u_{0}\left(  t\right)  =x_{0}^{0}\\%
\mathcal{L}%
\left[  u_{m+1}\left(  t\right)  \right]  =c_{0}\left[  \sum_{k=0}^{m}%
u_{k}^{\prime}\left(  t\right)  -f\left(  t,\sum_{k=0}^{m}u_{k}\left(
t\right)  \right)  \right]  ,m\geq0\\
\text{with }u_{k}\left(  t_{0}\right)  =0\text{, }\forall k\in%
\mathbb{N}%
\end{array}
\right. \label{eq_l_00500}%
\end{equation}
with $%
\mathcal{L}%
\left[  h\left(  t\right)  \right]  =\frac{dh}{dt}\left(  t\right)  $ and $%
\mathcal{L}%
^{-1}\left[  h\left(  t\right)  \right]  =\int_{t_{0}}^{t}h\left(  \xi \right)
d\xi$ and algorithm
\begin{equation}
\left \{
\begin{array}
[c]{l}%
x_{0}\left(  t\right)  =x_{0}^{\left(  0\right)  }\\
x_{m+1}=x_{0}^{\left(  0\right)  }+\left(  1+c_{0}\right)  \int_{t_{0}}%
^{t}\frac{dx_{m}}{dt}\left(  \xi \right)  d\xi-c_{0}\int_{t_{0}}^{t}f\left(
\xi,x_{m}\left(  \xi \right)  \right)  d\xi
\end{array}
\right.  .\label{eq_l_001000}%
\end{equation}
Then,
\begin{equation}
x_{m}\left(  t\right)  =\sum_{k=0}^{m}u_{k}\left(  t\right)  \text{, }\forall
m\in%
\mathbb{N}
_{0}\label{eq_l_002000}%
\end{equation}
whenever $\left(  t,x_{k}\left(  t\right)  \right)  \in D$ for $k=1,\ldots
,m-1.$

\begin{pf}
Let%
\'{}%
s argue by mathematical induction. For $m=0$ (\ref{eq_l_002000}) is trivially
true from definition. Consider now the inductive hypothesis. Suppose that
(\ref{eq_l_002000}) is true for some $p\in%
\mathbb{N}
$, that is, $x_{p}\left(  t\right)  =\sum_{k=0}^{p}u_{k}\left(  t\right)  $
and $\left(  t,x_{k}\left(  t\right)  \right)  \in D$ for $k=1,\ldots,p-1.$
Let's prove that
\[
u_{p+1}\left(  t\right)  =x_{p+1}\left(  t\right)  -x_{p}\left(  t\right)  ,
\]
that is, $x_{p+1}\left(  t\right)  =\sum_{k=0}^{p+1}u_{k}\left(  t\right)  $.
From (\ref{eq_l_00500}) and using the inductive hypothesis we successively
deduce%
\[
u_{p+1}^{\prime}\left(  t\right)  =c_{0}\left[  \sum_{k=0}^{p}u_{k}^{\prime
}\left(  t\right)  -f\left(  t,\sum_{k=0}^{p}u_{k}\left(  t\right)  \right)
\right]  \Rightarrow
\]%
\[
u_{p+1}^{\prime}\left(  t\right)  =c_{0}\left[  x_{p}^{\prime}\left(
t\right)  -f\left(  t,x_{p}\left(  t\right)  \right)  \right]  ,
\]
and from (\ref{eq_l_001000})%
\[
x_{p+1}^{\prime}\left(  t\right)  =\left(  1+c_{0}\right)  x_{p}^{\prime
}\left(  t\right)  -c_{0}f\left(  t,x_{p}\left(  t\right)  \right)
\Rightarrow
\]%
\[
\left(  x_{p+1}\left(  t\right)  -x_{p}\left(  t\right)  \right)  ^{\prime
}=c_{0} \left[  x_{p}^{\prime}\left(  t\right)  -f\left(  t,x_{p}\left(
t\right)  \right)  \right]  .
\]
Furthermore, $u_{p+1}\left(  t_{0}\right)  =0$ from (\ref{eq_l_00500}) and
$x_{p+1}\left(  t_{0}\right)  =x_{p}\left(  t_{0}\right)  =x_{0}^{\left(
0\right)  }$ from (\ref{eq_l_001000}), hence
\[
u_{p+1}\left(  t_{0}\right)  =x_{p+1}\left(  t_{0}\right)  -x_{p}\left(
t_{0}\right)  =0\text{.}%
\]
So, $u_{p+1}\left(  t\right)  =x_{p+1}\left(  t\right)  -x_{p}\left(
t\right)  $ for all $t$ such that $\left(  t,x_{p}\left(  t\right)  \right)
\in D$. This completes the inductive step.\hfill$\square$
\end{pf}
\end{thm}

We are interested in knowing for what values of $c_{0}$ can we guarantee the
convergence of the IFOHAM algorithm (\ref{eq_300010}) in the context of
choices (\ref{eq_300020}) and (\ref{eq_300030}). In this way, we will
establish some sufficient conditions for convergence of this algorithm.

Let us first present a trivial lemma that we will need.

\begin{lem}
\label{lema2}Let $\alpha$ and $\beta$ real constants and $h\left(  x\right)
=\left(  1+x\right)  \alpha-x\beta$ with $\left \vert \alpha \right \vert \leq A$
and $\left \vert \beta \right \vert \leq A.$ If $x\in \left[  -1,0\right]  $ then
$\left \vert h\left(  x\right)  \right \vert \leq A$.

\begin{pf}
Let $\alpha-\beta=\delta.$ Then, $h\left(  x\right)  =\alpha+x\delta$ and
$h\left(  x\right)  =\left(  1+x\right)  \delta+\beta.$ If $\delta=0$ then
$h\left(  x\right)  =\alpha=\beta$ $\forall x\in%
\mathbb{R}
$. Hence, $\left \vert h\left(  x\right)  \right \vert =\left \vert
\alpha \right \vert \leq A.$ If $\delta>0$ then $h\left(  x\right)
=\alpha+x\delta \leq \alpha$ and $h\left(  x\right)  =\left(  1+x\right)
\delta+\beta \geq \beta$ $\forall x\in \left[  -1,0\right]  $. Hence,
$-A\leq \beta \leq h\left(  x\right)  \leq \alpha \leq A$. Then, $\left \vert
h\left(  x\right)  \right \vert \leq A$. If $\delta<0$ then $h\left(  x\right)
=\alpha+x\delta \geq \alpha$ and $h\left(  x\right)  =\left(  1+x\right)
\delta+\beta \leq \beta$ $\forall x\in \left[  -1,0\right]  $. Hence,
$-A\leq \alpha \leq h\left(  x\right)  \leq \beta \leq A$. Then, $\left \vert
h\left(  x\right)  \right \vert \leq A$. So, If $x\in \left[  -1,0\right]  $
then $\left \vert h\left(  x\right)  \right \vert \leq A$.\hfill$\square$
\end{pf}
\end{lem}

\begin{thm}
\label{proposition_5}Let $D$ be an open set in $%
\mathbb{R}
^{2}$. Let $\left(  t_{0},x_{0}^{\left(  0\right)  }\right)  \in D$ and let
$a$ and $b$ be positive constants such that the set%
\[
R=\left \{  \left(  t,x\right)  :\left \vert t-t_{0}\right \vert \leq a\text{ and
}\left \vert x-x_{0}^{\left(  0\right)  }\right \vert \leq b\right \}
\]
is contained in $D.$ Suppose function $f$ is continuous and defined on $D$ and
satisfies a Lipschitz condition with respect to $x$ in $R$ with Lipschitz
constant $L$. Let $M=\underset{\left(  t,x\right)  \in R}{\max}\left \vert
f\left(  t,x\right)  \right \vert $ and $A=\min \left \{  a,\frac{b}{M}\right \}
.$ Consider the IVP
\begin{equation}
\left \{
\begin{array}
[c]{c}%
\frac{dx}{dt}=f\left(  t,x\right) \\
x\left(  t_{0}\right)  =x_{0}^{\left(  0\right)  }%
\end{array}
\right.  ,\label{eq_p4_0005}%
\end{equation}
and its unique solution $x=x\left(  t\right)  $ on the open interval
$I=\left]  t_{0}-A,t_{0}+A\right[  $. Consider also the IFOHAM algorithm
\begin{equation}
\left \{
\begin{array}
[c]{l}%
u_{0}\left(  t\right)  =x_{0}^{0}\\%
\mathcal{L}%
\left[  u_{m+1}\left(  t\right)  \right]  =c_{0}\left[  \sum_{k=0}^{m}%
u_{k}^{\prime}\left(  t\right)  -f\left(  t,\sum_{k=0}^{m}u_{k}\left(
t\right)  \right)  \right]  ,m\geq0\\
\text{with }u_{k}\left(  t_{0}\right)  =0\text{, }\forall k\in%
\mathbb{N}%
\end{array}
\right.  ,\label{eq_p4_0010}%
\end{equation}
with $%
\mathcal{L}%
\left[  h\left(  t\right)  \right]  =\frac{dh}{dt}\left(  t\right)  $ and $%
\mathcal{L}%
^{-1}\left[  h\left(  t\right)  \right]  =\int_{t_{0}}^{t}h\left(  \xi \right)
d\xi$ and algorithm
\begin{equation}
\left \{
\begin{array}
[c]{l}%
x_{0}\left(  t\right)  =x_{0}^{\left(  0\right)  }\\
x_{m+1}=x_{0}^{\left(  0\right)  }+\left(  1+c_{0}\right)  \int_{t_{0}}%
^{t}\frac{dx_{m}}{dt}\left(  \xi \right)  d\xi-c_{0}\int_{t_{0}}^{t}f\left(
\xi,x_{m}\left(  \xi \right)  \right)  d\xi
\end{array}
\right. \label{eq_p4_0015}%
\end{equation}
and its associated operator
\begin{equation}
F\left(  x\left(  t\right)  \right)  =\left(  1+c_{0}\right)  x\left(
t\right)  -c_{0}\left(  x_{0}^{\left(  0\right)  }+\int_{t_{0}}^{t}f\left(
\xi,x\left(  \xi \right)  \right)  d\xi \right)  .\label{eq_p5_001010}%
\end{equation}

\begin{enumerate}
\item If $c_{0}\in \left[  -1,0\right[  $ then $\left \{  \sum_{k=0}^{n}%
u_{k}\left(  t\right)  \right \}  $ converges uniformly on $I$ to $x\left(
t\right)  $.

\item Define%
\[
S=\left \{  x\left(  t\right)  \in C\left(  J\right)  :\left \vert x\left(
t\right)  -x_{0}^{\left(  0\right)  }\right \vert \leq b\text{, }\left \vert
t-t_{0}\right \vert \leq A\right \}  ,
\]
let $\tilde{L}$ be any constant $\tilde{L}>L,$ $J=\left[  t_{0}-A,t_{0}%
+A\right]  $ and consider the norm defined as follows:%
\[
\left \Vert x\right \Vert _{e}=\max_{t\in J}\left \vert x\left(  t\right)
e^{-\tilde{L}\left \vert t-t_{0}\right \vert }\right \vert .
\]

If $x\left(  t\right)  $ and $y\left(  t\right)  $ belongs to $S$ and
$c_{0}\in \left[  -1,0\right[  $ then
\[
\left \Vert F\left(  x\left(  t\right)  \right)  -F\left(  y\left(  t\right)
\right)  \right \Vert _{e}\leq k\left \Vert x\left(  t\right)  -y\left(
t\right)  \right \Vert _{e}%
\]
with
\begin{equation}
0<k=1+\left(  1-\frac{L}{\tilde{L}}\left(  1-e^{-\tilde{L}A}\right)  \right)
c_{0}<1.\label{eq_p5_002020}%
\end{equation}

\end{enumerate}

\begin{pf}
We begin by demonstrating part 1. The demonstration of part 2 will follow from
the latter. From Proposition \ref{proposition_4} one knows that%
\[
x_{n}\left(  t\right)  =\sum_{k=0}^{n}u_{k}\left(  t\right)  .
\]
Hence, it is sufficient to show that if $c_{0}\in \left[  -1,0\right[  $ then
$\left \{  x_{n}\left(  t\right)  \right \}  $ converges uniformly on $I$ to
$x\left(  t\right)  $. Therefore consider algorithm (\ref{eq_p4_0015})and its
associated operator (\ref{eq_p5_001010}).

Let's show that $\left \{  x_{n}\left(  t\right)  \right \}  $ converges
uniformly on $I$ to some $y\left(  t\right)  $ using Banach's fixed point
theorem. The missing details of this elementary approach can be found in
\cite{kreyszig1} and \cite{zeidler1}, for instance.

Let $J=\left[  t_{0}-A,t_{0}+A\right]  $ and define the (non empty, closed)
subset%
\[
S=\left \{  x\in C\left(  J\right)  :\left \vert x\left(  t\right)
-x_{0}^{\left(  0\right)  }\right \vert \leq b\text{, }\left \vert
t-t_{0}\right \vert \leq A\right \}  ,
\]

of the Banach space $C\left(  J\right)  $ with the norm $\left \Vert
\cdot \right \Vert _{\infty}$. Note that
\[
\left \vert F\left(  x\left(  t\right)  \right)  -x_{0}^{\left(  0\right)
}\right \vert =\left \vert \left(  1+c_{0}\right)  \left(  x\left(  t\right)
-x_{0}^{\left(  0\right)  }\right)  -c_{0}\int_{t_{0}}^{t}f\left(
\xi,x\left(  \xi \right)  \right)  d\xi \right \vert .
\]
If $x\left(  t\right)  \in S,$ then $\left \Vert x-x_{0}^{\left(  0\right)
}\right \Vert _{\infty}\leq b$ and
\[
\left \Vert \int_{t_{0}}^{t}f\left(  \xi,x\left(  \xi \right)  \right)
d\xi \right \Vert _{\infty}=\max_{t\in J}\left \vert \int_{t_{0}}^{t}f\left(
\xi,x\left(  \xi \right)  \right)  d\xi \right \vert \leq AM\leq b\text{. }%
\]
Note also that $c_{0}\in \left[  -1,0\right[  ,$ so we can conclude from Lemma
(\ref{lema2}) that $$\left \Vert F\left(  x\left(  t\right)  \right)
-x_{0}^{\left(  0\right)  }\right \Vert _{\infty}\leq b.$$ Hence,
$$x\left(  t\right)  \in S\Rightarrow F\left(  x\left(  t\right)  \right)  \in
S.$$

Let $\tilde{L}$ be any constant $\tilde{L}>L$ and consider the norm %
\[
\left \Vert x\right \Vert _{e}=\max_{t\in J}\left \vert x\left(  t\right)
e^{-\tilde{L}\left \vert t-t_{0}\right \vert }\right \vert .
\]%
Observe that norms $\left \Vert \cdot \right \Vert _{e}$ and $\left \Vert
\cdot \right \Vert _{\infty}$ are equivalent. Suppose $x\left(  t\right)  $ and
$y\left(  t\right)  $ are in $S$ and consider now the expression%
\begin{gather*}
\left \vert F\left(  x\left(  t\right)  \right)  -F\left(  y\left(  t\right)
\right)  \right \vert =\\
=\left \vert \left(  1+c_{0}\right)  \left(  x\left(  t\right)  -y\left(
t\right)  \right)  -c_{0}\left(  \int_{t_{0}}^{t}\left(  f\left(  \xi,x\left(
\xi \right)  \right)  -f\left(  \xi,y\left(  \xi \right)  \right)  \right)
d\xi \right)  \right \vert ,
\end{gather*}
obtained from (\ref{eq_p5_001010}). Clearly%
\begin{gather*}
\left \vert F\left(  x\left(  t\right)  \right)  -F\left(  y\left(  t\right)
\right)  \right \vert e^{-\tilde{L}\left \vert t-t_{0}\right \vert }\leq \\
\leq \left \vert 1+c_{0}\right \vert \left \vert x\left(  t\right)  -y\left(
t\right)  \right \vert e^{-\tilde{L}\left \vert t-t_{0}\right \vert }+\left \vert
c_{0}\right \vert L\left \vert \int_{t_{0}}^{t}\left \vert x\left(  \xi \right)
-y\left(  \xi \right)  \right \vert d\xi \right \vert e^{-\tilde{L}\left \vert
t-t_{0}\right \vert }%
\end{gather*}
and%
\begin{gather*}
\left \Vert F\left(  x\left(  t\right)  \right)  -F\left(  y\left(  t\right)
\right)  \right \Vert _{e}\leq \\
\leq \left \vert 1+c_{0}\right \vert \left \Vert x\left(  t\right)  -y\left(
t\right)  \right \Vert _{e}+\left \vert c_{0}\right \vert L\left \Vert \int
_{t_{0}}^{t}\left \vert x\left(  \xi \right)  -y\left(  \xi \right)  \right \vert
d\xi \right \Vert _{e}.
\end{gather*}
One can deduce that%
\[
\left \Vert \int_{t_{0}}^{t}\left \vert x\left(  \xi \right)  -y\left(
\xi \right)  \right \vert d\xi \right \Vert _{e}\leq \frac{\left \Vert x\left(
t\right)  -y\left(  t\right)  \right \Vert _{e}}{\tilde{L}}\left(
1-e^{-\tilde{L}A}\right)  ,
\]
therefore%
\begin{gather*}
\left \Vert F\left(  x\left(  t\right)  \right)  -F\left(  y\left(  t\right)
\right)  \right \Vert _{e}\leq \\
\leq \left \{  \left \vert 1+c_{0}\right \vert +\left \vert c_{0}\right \vert
\frac{L}{\tilde{L}}\left(  1-e^{-\tilde{L}A}\right)  \right \}  \left \Vert
x\left(  t\right)  -y\left(  t\right)  \right \Vert _{e}.
\end{gather*}
If $c_{0}\in \left[  -1,0\right[  $, observe that
\[
0<k=\left \vert 1+c_{0}\right \vert +\left \vert c_{0}\right \vert \frac{L}%
{\tilde{L}}\left(  1-e^{-\tilde{L}A}\right)  =\left(  1+c_{0}\right)
-c_{0}\frac{L}{\tilde{L}}\left(  1-e^{-\tilde{L}A}\right)  <1,
\]
that is
\[
0<1+\left(  1-\frac{L}{\tilde{L}}\left(  1-e^{-\tilde{L}A}\right)  \right)
c_{0}<1.
\]
So,
\[
\left \Vert F\left(  x\left(  t\right)  \right)  -F\left(  y\left(  t\right)
\right)  \right \Vert _{e}\leq k\left \Vert x\left(  t\right)  -y\left(
t\right)  \right \Vert _{e},
\]
is a contraction. Therefore, from Banach's fixed point theorem one conclude
that $\left \{  x_{n}\left(  t\right)  \right \}  $ converges uniformly on $J$
to some fixed point $y\left(  t\right)  $ of (\ref{eq_p5_001010}). Clearly, if
$y\left(  t\right)  $ is the fixed point of (\ref{eq_p5_001010}) then, one
deduce also that
\[
y\left(  t\right)  =x_{0}^{\left(  0\right)  }+\int_{t_{0}}^{t}f\left(
\xi,y\left(  \xi \right)  \right)  d\xi,
\]
that is, $y\left(  t\right)  $ is the solution the IVP (\ref{eq_p4_0005}) on
the interior of $J$. From the uniqueness of the solution we will conclude that
$y\left(  t\right)  =x\left(  t\right)  $ on $I$. This completes the proof of
both parts.\hfill$\square$
\end{pf}
\end{thm}

We would like to stress that Proposition \ref{proposition_5} establishes
sufficient conditions for the convergence of IFOHAM under the corresponding
context. The convergence also depends on the structure of $f$. So, it will not
come as a surprise if convergence is also verified in a wider range $\left[
c,0\right[  $ with $c<-1$.

Moreover, expression (\ref{eq_p5_002020}) suggest that the minimum on $\left[
-1,0\right[  $ of the contraction constant $k$ is attained at $c_{0}=-1$. This
means that in this frame and in the absence of information about the
convergence of IFOHAM for $c_{0}$ less than $-1$ the best choice for this
parameter will be $c_{0}=-1$, that is, the best choice will be
Picard-Lindel\"{o}ff 's iteration algorithm. So, the knowledge of the
structure of $f$ in (\ref{eq_300005}) is of primordial importance for the
useful use of the IFOHAM algorithm in the studied context.

\section{Results and discussion}

In order to preliminary compare the relative performance of HAM and IFOHAM we
will address again the IVP (\ref{eq_001010}).

\begin{figure}
[t]
\begin{center}
\includegraphics[totalheight=4in
]%
{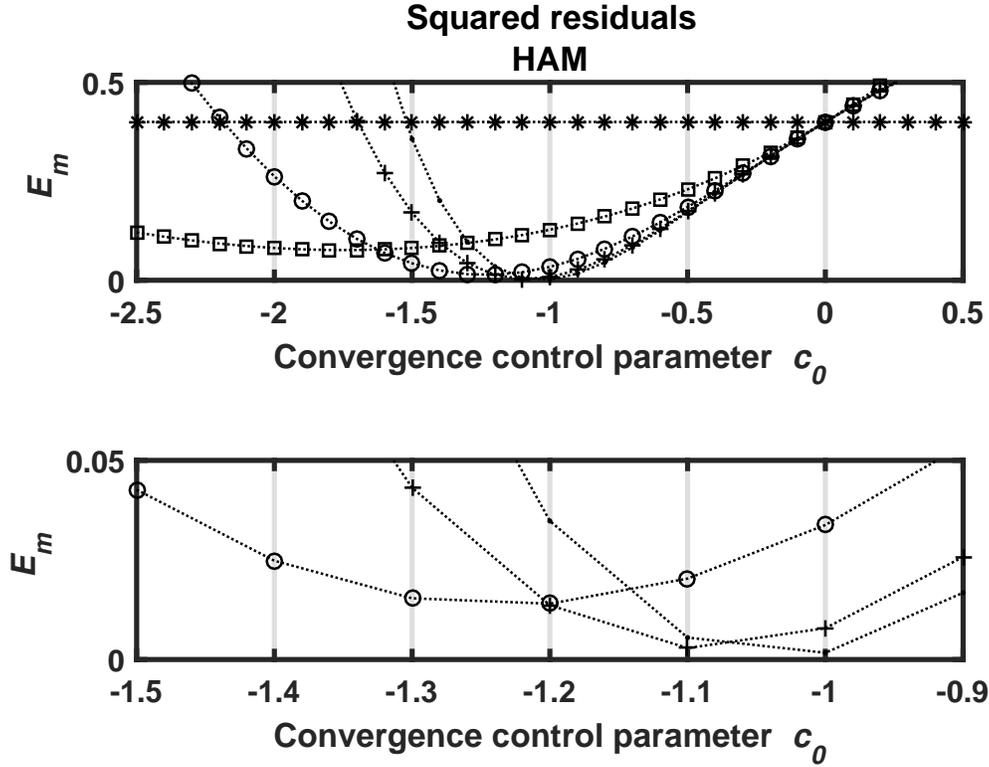}%
\caption{Approximate solutions: $\ast$-zeroth order, $\square$- first order,
$\bigcirc$-second order, $+$ - third order, $\bullet$ - fourth order.}%
\label{1fig_dupla_ham_conv}%
\end{center}
\end{figure}

\begin{figure}
[t]
\begin{center}
\includegraphics[totalheight=4in
]%
{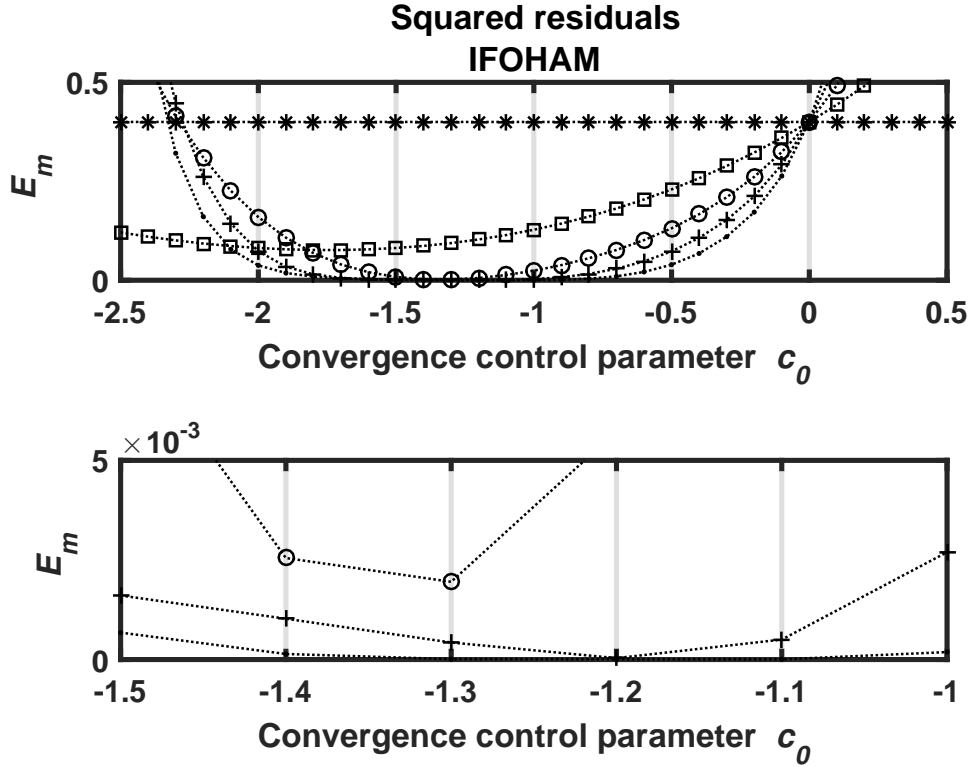}%
\caption{Approximate solutions: $\ast$-zeroth order, $\square$- first order,
$\bigcirc$-second order, $+$ - third order, $\bullet$ - fourth order.}%
\label{1fig_dupla_ifoham_conv}%
\end{center}
\end{figure}

In Figures \ref{1fig_dupla_ham_conv} and \ref{1fig_dupla_ifoham_conv} we
display for different values of the convergence control parameter $c_{0}$ the
squared residuals $E_{m}$ corresponding to different $m$th-order solutions
obtained using HAM and IFOHAM. The squared residuals were computed using
expression%
\[
E_{m}=\int_{\Omega}N\left(  \sum_{i=0}^{m}u_{i}\right)  dt
\]
where $N$ represents operator (\ref{eq_001005} ) and $\Omega=\left[
-1,1\right]  $. In the bottom sub-figures we display a more detailed zoom to
improve the determination of the location of the value of the parameter
$c_{0}$ that minimizes $E_{m}$.%

With respect to Figure \ref{1fig_dupla_ham_conv} and concerning the HAM, data
suggest that:

\begin{itemize}
\item HAM converges for $c_{0}\in \left[  -1\text{,}0\right[  $ and diverges
por $c_{0}>0$;

\item Performance of the HAM\ algorithm for this test case improves in the
neighborhood of $c_{0}=-1.$
\end{itemize}%

With respect to Figure \ref{1fig_dupla_ifoham_conv} and concerning IFOHAM,
data suggest that:

\begin{itemize}
\item IFOHAM converges for $c_{0}\in \left[  -1.3,\text{ }0\right[  $ and
diverges por $c_{0}\geq0$;

\item Performance of the IFOHAM\ algorithm is the best in the neighborhood
$c_{0}=-1.2.$
\end{itemize}

Note that the convergence of IFOHAM is assured if $c_{0}\in \left[
-1\text{,}0\right[  $ in agreement with Proposition \ref{proposition_5}.
However, depending on the struture of $f$ in (\ref{eq_300005}), convergence of
IFOHAM, as noted in this case, can occur over a wider range $\left[
c,0\right[  $ with $c<-1$.

One observe also that, the performance of IFOHAM, for $c_{0}\in \left[
-1,\text{ }0\right[  ,$ is best at the left end of this range. This fact is in
agreement with expression (\ref{eq_p5_002020}) since the minimum value of the
contraction constant $k$ on $\left[  -1,\text{ }0\right[  $ interval is
attained at $c_{0}=-1$. As previously mentioned at the end of the last
section, this means that in the absence of information about the convergence
of IFOHAM for $c_{0}$ less than $-1,$ the best choice for this parameter will
be $c_{0}=-1$, that is, the best choice will be Picard-Lindel\"{o}ff 's
iteration algorithm. So, the knowledge of the structure of $f$ in
(\ref{eq_300005}) is essential for an effective use of the IFOHAM algorithm in
the studied context.

In Tables \ref{HAM_table1}, \ref{IFOHAM_table2} and \ref{IFOHAM_table3} we
display the computed squared residuals $E_{m}$ as well as the computational
CPU time consumed to obtain the corresponding $m$th-order approximate
solutions for cases $c_{0}=-1$ (HAM), $c_{0}=-1$ (IFOHAM) and $c_{0}=-1.2$
(IFOHAM). The above cases have been chosen especially because:

\begin{itemize}
\item HAM is better effective in the neighborhood of $c_{0}=-1$ as was
suggested from the analysis of Figure \ref{1fig_dupla_ham_conv};

\item IFOHAM with $c_{0}=-1$ (that is, Picard-Lindel\"{o}ff 's iteration
algorithm) is the best blind implementation of IFOHAM in the absence of
information regarding the structure of $f$;

\item IFOHAM in the neighborhood of $c_{0}=-1.2$ is the best informed
implementation of IFOHAM as was suggested from the analysis of Figure
\ref{1fig_dupla_ifoham_conv}.
\end{itemize}

Considering the extension of\ some expressions of the $m$th-order terms and
$m$th-order approximate solutions these expressions were only partially
reproduced in the Tables \ref{IFOHAM_table2} and \ref{IFOHAM_table3}. However,
the missing terms replaced by suspension points can be easily obtained by
applying the IFOHAM technique on a symbolic computer environment.%

\begin{table}[tbp] \centering
\caption{HAM effectiveness, $c_{0}$ =  $-1$ in computing the $m$th-order
approximate solution}%
\begin{tabular}
[c]{lrccl}\hline \hline
$k/m$ & $u_{k}$ & $\sum_{k=0}^{m}u_{k}$ & $E_{m}$ & CPU time [s]\\ \hline \hline
\multicolumn{1}{c}{$0$} & $t$ & \multicolumn{1}{r}{$t$} & $4.00e-01$ &
\multicolumn{1}{c}{$0.000$}\\
\multicolumn{1}{c}{$1$} & $\frac{t^{3}}{3}$ & \multicolumn{1}{r}{$t+\frac
{t^{3}}{3}$} & $1.28e-01$ & \multicolumn{1}{c}{$2.109$}\\
\multicolumn{1}{c}{$2$} & $\frac{2t^{5}}{15}$ & \multicolumn{1}{r}{$t+\frac
{t^{3}}{3}+\frac{2t^{5}}{15}$} & $3.38e-02$ & \multicolumn{1}{c}{$2.422$}\\
\multicolumn{1}{c}{$3$} & $\frac{17t^{7}}{315}$ & \multicolumn{1}{r}{$t+\frac
{t^{3}}{3}+\frac{2t^{5}}{15}+\frac{17t^{7}}{315}$} & $7.88e-03$ &
\multicolumn{1}{c}{$2.719$}\\
\multicolumn{1}{c}{$4$} & $\frac{62t^{9}}{2835}$ & \multicolumn{1}{r}{$t+\frac
{t^{3}}{3}+\frac{2t^{5}}{15}+\frac{17t^{7}}{315}+\frac{62t^{9}}{2835}$} &
$1.70e-03$ & \multicolumn{1}{c}{$2.984$}\\ \hline \hline
\end{tabular}
\label{HAM_table1}%
\end{table}%

The tabulated data suggest that in addressing our test case, the IVP (\ref{eq_001010}),
Picard-Lindel\"{o}ff 's iteration algorithm (IFOHAM with $c_{0}=-1$) is better
effective than the best implementation of HAM (HAM with $c_{0}=-1$) and the
implementation of IFOHAM\ with $c_{0}=-1.2$ is the best of all the illustrated implementations.%

\begin{table}[tbp] \centering
\caption{IFOHAM effectiveness, $c_{0}$ =  $-1$ in computing the $m$th-order
approximate solution}%
\begin{tabular}
[c]{crrr}\hline \hline
$k$ & $u_{k}$ &  & \\ \hline \hline
$0$ & $t$ &  & \\
$1$ & $\frac{t^{3}}{3}$ &  & \\
$2$ & $\frac{t^{7}}{63}+\frac{2t^{5}}{15}$ &  & \\
$3$ & $\frac{t^{15}}{59535}+\frac{4t^{13}}{12285}+\frac{134t^{11}}%
{51975}+\frac{38t^{9}}{2835}+\frac{4t^{7}}{105}$ &  & \\
$4$ & $\frac{t^{31}}{109876902975}+\cdots+\frac{8t^{9}}{945}$ &  &
\\ \hline \hline
$m$ & \multicolumn{1}{c}{$\sum_{k=0}^{m}u_{k}$} & \multicolumn{1}{c}{$E_{m}$}
& \multicolumn{1}{c}{CPU time [s]}\\ \hline \hline
$0$ & $t$ & \multicolumn{1}{c}{$4.00e-01$} & \multicolumn{1}{c}{$0.000$}\\
$1$ & $\frac{t^{3}}{3}+t$ & \multicolumn{1}{c}{$1.28e-01$} &
\multicolumn{1}{c}{$0.609$}\\
$2$ & $\frac{t^{7}}{63}+\frac{2t^{5}}{15}+\frac{t^{3}}{3}+t$ &
\multicolumn{1}{c}{$2.42e-02$} & \multicolumn{1}{c}{$1.031$}\\
$3$ & $\frac{t^{15}}{59535}+\frac{4t^{13}}{12285}+\cdots+\frac{17}{315}%
t^{7}+\frac{2t^{5}}{15}+\frac{t^{3}}{3}+t$ & \multicolumn{1}{c}{$2.69e-03$} &
\multicolumn{1}{c}{$1.406$}\\
$4$ & $\frac{t^{31}}{109876902\,975}+\cdots+\frac{t^{3}}{3}+t$ &
\multicolumn{1}{c}{$1.87e-04$} & \multicolumn{1}{c}{$1.938$}\\ \hline \hline
\end{tabular}
\label{IFOHAM_table2}%
\end{table}%
%

\begin{table}[tbp] \centering
\caption{IFOHAM effectiveness, $c_{0}$ =  $-1.2$ in computing the
$m$th-order approximate solution}%
\begin{tabular}
[c]{crrr}\hline \hline
$k$ & $u_{k}$ &  & \\ \hline \hline
$0$ & $t$ &  & \\
$1$ & $\frac{2t^{3}}{5}$ &  & \\
$2$ & $\frac{24t^{7}}{875}+\frac{24t^{5}}{125}-\frac{2t^{3}}{25}$ &  & \\
$3$ & $\frac{1152t^{15}}{19140625}+\cdots+\frac{1104t^{7}}{21875}%
-\frac{48t^{5}}{625}+\frac{2t^{3}}{125}$ &  & \\
$4$ & $\frac{7962624t^{31}}{56786346435546875}+\cdots+\frac{72t^{5}}%
{3125}-\frac{2t^{3}}{625}$ &  & \\ \hline \hline
$M$ & \multicolumn{1}{c}{$\sum_{k=0}^{M}u_{k}$} & \multicolumn{1}{c}{$E_{M}$}
& \multicolumn{1}{c}{CPU time [s]}\\ \hline \hline
$0$ & $t$ & \multicolumn{1}{c}{$4.00e-01$} & \multicolumn{1}{c}{$0.000$}\\
$1$ & $\frac{2t^{3}}{5}+t$ & \multicolumn{1}{c}{$1.03e-01$} &
\multicolumn{1}{c}{$0.609$}\\
$2$ & $\frac{24t^{7}}{875}+\frac{24t^{5}}{125}+\frac{8t^{3}}{25}+t$ &
\multicolumn{1}{c}{$5.73e-03$} & \multicolumn{1}{c}{$1.063$}\\
$3$ & $\frac{1152t^{15}}{19\,140\,625}+\cdots+\frac{1704t^{7}}{21\,875}%
+\frac{72t^{5}}{625}+\frac{42t^{3}}{125}+t$ & \multicolumn{1}{c}{$3.54e-05$} &
\multicolumn{1}{c}{$1.438$}\\
$4$ & $\frac{7962624t^{31}}{56786346435546875}+\cdots+\frac{432t^{5}}%
{3125}+\frac{208t^{3}}{625}+t$ & \multicolumn{1}{c}{$5.45e-06$} &
\multicolumn{1}{c}{$2.078$}\\ \hline \hline
\end{tabular}
\label{IFOHAM_table3}%
\end{table}%

Note that sequences of approximate solutions generated by HAM or IFOHAM
converge to the MacLaurin series of $x=\tan t$ (the exact known solution of
our problem). Despite this fact, it should be noted that the terms of each
approximate solution already calculated in one iteration using IFOHAM may be
modified in the next iteration contrary to what happens using HAM. As was
noted before, HAM can \textquotedblleft surgically\textquotedblright \ determines %
the terms of the Maclaurin series of the solution of our problem.

Moreover, in a few iterations the IFOHAM algorithm has to handle particularly
long expressions. This may constitute a drawback of this algorithm.

However, these preliminary tests suggest that IFOHAM exhibits an interesting
performance both in aspects related to the speed of convergence and in aspects
related to the CPU calculation time.

\section{Conclusion and future work}

In addressing the classic IVP problem%
\begin{equation}
\left \{
\begin{array}
[c]{c}%
\frac{dx}{dt}=f\left(  t,x\right) \\
x\left(  t_{0}\right)  =x_{0}^{\left(  0\right)  }%
\end{array}
\right.  ,
\end{equation}
we found that, conveniently defining $%
\mathcal{L}%
\left[  h\left(  t\right)  \right]  =\frac{dh}{dt}\left(  t\right)  ,$ IFOHAM%
\begin{equation}%
\mathcal{L}%
\left[  u_{m+1}\left(  t\right)  \right]  =c_{0}\left[  N\left[
\sum_{k=0}^{m}u_{k}\left(  t\right)  \right]  \right]  \text{, }m\geq0,
\end{equation}
with $c_{0}=-1$ coincides exactly with Picard-Lindel\"{o}ff 's iteration
algorithm. We concluded also that IFOHAM converges if $c_{0}\in \left]
-1,0\right[  $ and depending on the structure of $f$ IFOHAM can still converge
with a better convergence speed to the searched solution if $c_{0}<-1$.
Clearly, the knowledge of the structure of $f$ is of primordial importance for
the future useful use of the IFOHAM algorithm in the studied context. Given
these facts one can state that IFOHAM generalizes Picard-Lindel\"{o}ff 's
iteration algorithm.

Preliminary tests showed that IFOHAM exhibited a very good performance both in
aspects related to the speed of convergence and in aspects related to the CPU
calculation time.

A very favorable aspect of IFOHAM lies in the ease of its implementation which
is simple and without complexities. However, in a few iterations the IFOHAM
algorithm has to handle particularly long expressions. This may constitute a
drawback of this algorithm.

With regard to future work we would like to mention some possible interesting
directions we are presently dealing with:

\begin{itemize}
\item To study the convergence of IFOHAM with respect the structure of $f$ in
(\ref{eq_300005}) or more generally regarding the structure of the operator
$N$ in (\ref{eq_0010});

\item To study the existence of flexibility of IFOHAM on the choice of base
functions and decide about the solution expression by adequately choosing $%
\mathcal{L}%
$ and the initial guess $u_{0}\left(  t\right)  $ as in the use of HAM;

\item To study the ability of IFOHAM to find the main parameters, such as
amplitude and frequency, of periodic solutions of nonlinear evolution problems;

\item Study of the applicability of IFOHAM in addressing other classes of evolution non-linear problems.
\end{itemize}

\section{Acknowledgments}

We would like to express our acknowledgments to my colleague Professor
M\'{a}rio Gatta by the interesting discussions concerning this work.


\end{document}